\documentclass[a4paper,twoside,11pt]{article}
\usepackage{a4,color,graphics,lmodern,fancyhdr}

\usepackage{amsfonts}
\def\BbR{{\mathbb R}}

\usepackage[latin1]{inputenc}
\usepackage[T1]{fontenc}
\usepackage[francais]{babel}
\usepackage{amsmath}
\usepackage{amsthm}
\usepackage{amsfonts}
\usepackage{hyperref}
\begin{document}
\begin{titlepage}
%% Title Page %%%%%%%%%%%%%%%%%%%%%%%%%%%%%%%%%%%%%%%%%%%%%%%
%% ==> Write your text here or include other files.
%\Large
\begin{center}
\Large
%\textsc\textbf{MINISTERE DU TRANSPORT ET DE L'EQUIPEMENT}\\
%\%textsc{Office de la Topographie et du Cadastre }\\

\vspace{6cm}

\LARGE
\textsc{\textbf{La G\'eom\'etrie de Compensation Non-Lin\'eaire - Le Probl\`eme Spatial d'Intersection \\  dans l'Option de la G\'eod\'esie Tridimensionnelle}}\\
[0.5\baselineskip]
%de\\[0.5\baselineskip]
%\\[0.5\baselineskip]
{Par \\ Abdelmajid BEN HADJ SALEM}\\
\vspace{0.5cm}
\normalsize
\textsc{Ing\'enieur G\'en\'eral  Retrait\'e de  l'Office de la Topographie et du Cadastre}\\

\vspace{1cm}
\textsc{Janvier 2017}
\\ %%Date - better you write it yourself.

\vspace{1cm}
\textsc{Version 1.}\\
%Prof. V. Nachname\\
%Dr. V. Nachname}\\

\vspace{5cm}
\textsc{abenhadjsalem@gmail.com}
%Faculty of ...\\
%Institute of ...}\\

\end{center}
%\vspace{2cm}
%% The simple version:
\end{titlepage}
%\huge
%\part{GEODESIE CLASSIQUE}
\tableofcontents
\clearpage %Table of contents
\Large
\begin{center}
\textbf{La G\'eom\'etrie de Compensation Non-Lin\'eaire - Le Probl\`eme Spatial d'Intersection \\ dans l'Option de la G\'eod\'esie Tridimensionnelle}
\\
\vspace{2mm}

\textsc{\textbf{Abdelmajid Ben Hadj Salem, Dipl. Ing.}\footnote{6, rue du Nil, Cit\'e Soliman Er-Riadh, 8020 Soliman, Tunisie }}

\end{center}
\date{}
%\footnotetext[1]{Partially supported by the Russian Basic Research Foundation,
%project {\bf 96-01-01104} and Institut Universitaire de France}

%\maketitle

%\centerline{to S.P. Novikov admiringly}
\normalsize
\vspace{6mm}
\section{Introduction}

%\vspace{4mm}
Dans un article ~\cite{bib1} E. Grafarend et B. Schaffrin ont \'etudi\'e la g\'eom\'etrie de la compensation ou ajustement non-lin\'eaire pour le cas du probl\`eme d'intersection plane en utilisant le mod\`ele de Gauss Markov, par les moindres carr\'es. Le pr\'esent papier d\'eveloppe la m\^eme m\'ethode en travaillant sur un exemple de la d\'etermination d'un point par trilat\'eration dans l'option de la g\'eod\'esie tridimensionnelle pour la d\'etermination des coordonn\'ees $(x,y,z)$ d'un point inconnu \`a partir des mesures des distances vers $n$ points connus.

\vspace{4mm}
\section{La G\'eom\'etrie Non Lin\'eaire du Mod\`ele de Gauss-Markov}
Le mod\`ele non lin\'eaire de Gauss-Markov est d\'efini par:
\begin{equation}
	\zeta(X)=L-e; \quad e \in \mathcal{N}(0,\Gamma) \label{m1}
\end{equation}
avec:

- $L$: le vecteur des observations $(n\times1)=(L_1,L_2,..,L_n)^T$,

- $X$: le vecteur des inconnues $(m\times1)=(X_1,X_2,..,X_m)^T$,

- $e$: le vecteur des erreurs $(n\times1)=(e_1,e_2,..,e_n)^T$ suit la loi normale $\mathcal{N}(0,\Gamma)$ avec $E(e)=0$ et $\Gamma=E(ee^T)$ la matrice de dispersion ou variance, on prendra $\Gamma=\sigma^2_0.P^{-1}$. $P$ est la matrice des poids et $\sigma_0$ une constante positive.

- $\zeta$: est une fonction donn\'ee injective d'un ouvert $U\subset \BbR^m \rightarrow \BbR^n$ et  $m<n$.
\\

Remarque: dans le cas d'un mod\`ele lin\'eaire, la fonction $\zeta=A.X$ o\`u $A$ est une matrice $n\times m$.
\\

On note $Im \zeta=\left\{\zeta(X) / X \in U\right\}$ l'image de $U$ par la fonction $\zeta$. $Im\zeta$ est une vari\'et\'e de dimension $m$ v\'erifiant les conditions:
\vspace{4mm}

(i): les vecteurs $\displaystyle \frac{\partial \zeta}{\partial X_1},\frac{\partial \zeta}{\partial X_2},...,\frac{\partial \zeta}{\partial X_m}$ sont lin\'eairement ind\'ependants en chaque point $X \in U$,
\vspace{4mm}

(ii): les fonctions $\displaystyle \frac{\partial^2 \zeta}{\partial X_i \partial X_j}$ sont continues sur $U$ pour $i,j \in \left\{1,2,...,m\right\}$.
\\

On introduit un produit scalaire:
\begin{equation}
	<\zeta_1,\zeta_2>=\zeta_1^T.P.\zeta_2 \label{m2}
\end{equation}
 D'o\`u la norme du vecteur $\zeta=(\zeta_1,\zeta_2,...,\zeta_n)^T$:
 \begin{equation}
	\|\zeta \|^2=<\zeta,\zeta>=\zeta^T.P.\zeta =\sum_{i=1}^np_i.\zeta_i^2 \label{m3}
\end{equation}
dans l'espace vectoriel $\BbR^n$ en prenant la matrice de poids $P$ une matrice diagonale . 
\\

Alors la solution par les moindres carr\'es $\bar{X}$ sera d\'efinie par:
\begin{equation}
		\|L-\bar{\zeta}(\bar{X}) \|=min\left\{\|L-\zeta(X) \| \,/ \,X\in U \right\} \label{m4}
\end{equation}
Cette condition est exprim\'ee par les \'equations de Lagrange-Euler soit:
\begin{equation}
	\frac{\partial}{\partial X_i} \|L-\zeta(X) \|^2=0 \quad \mbox{pour}\,i\in \left\{1,2,...,m\right\} \label{m5}
\end{equation}
En effet, on veut minimiser la fonction:
\begin{equation}
	F(X)=F(X_1,X_2,..,X_m)=\|L-\zeta(X) \|=\|L-\zeta(X_1,X_2,...,X_m) \|
\end{equation}
Comme $F$ est une fonction positive, minimiser $F$ c'est aussi minimiser $F^2$, soit $G(X)=F^2(X)$. En appliquant les \'equations de Lagrange-Euler, on obtient:
$$	-\frac{\partial G(X)}{\partial X_i}=0 \Rightarrow \frac{\partial G(X)}{\partial X_i}=0 $$ 
soit:
\begin{equation}
	\frac{\partial}{\partial X_i} \|L-\zeta(X_1,X_2,...,X_m) \|^2=0 \quad \mbox{pour}\,i\in \left\{1,2,...,m\right\} \label{m6}
\end{equation}
or:
\begin{align}
	 \|L-\zeta(X_1,X_2,...,X_m) \|^2=(L-\zeta(X_1,X_2,...,X_m))^T.P.(L-\zeta(X_1,X_2,...,X_m))=\nonumber \\
	 \zeta(X)^T.P.\zeta(X) -2L^T.P.\zeta(X)+L^T.P.L \label{m7} 
	\end{align}
Soit :
\begin{equation}
		\frac{\partial  G(X)}{\partial X_i}=2\zeta(X)^T.P.	\frac{\partial \zeta(X)}{\partial X_i}-2L^T.P.\frac{\partial \zeta(X)}{\partial X_i} \quad \mbox{pour}\,i\in \left\{1,2,...,m\right\} \label{m7a}
\end{equation}
ou encore :
$$ 		\frac{\partial  G(X)}{\partial X_i}=2(\zeta(X)-L)^T.P.\frac{\partial \zeta(X)}{\partial X_i} \quad \mbox{pour}\,i\in \left\{1,2,...,m\right\} $$

ce qui donne en utilisant (\ref{m5}):
\begin{align}
<L-\zeta(X),\frac{\partial \zeta(X)}{\partial X_i}>=0 \quad \mbox{pour}\,i\in \left\{1,2,...,m\right\} \label{m8} \\
\mbox{ou}\quad <e,\frac{\partial \zeta(X)}{\partial X_i}>=0 \quad \mbox{pour}\,i\in \left\{1,2,...,m\right\} \label{m9}
\end{align}
G\'eom\'etriquement, cela veut dire que le vecteur erreur $e =L-\zeta(X)$ est perpendiculaire (produit scalaire nul) au plan tangent de la vari\'et\'e $Im \zeta$ au point $\bar{\zeta}(\bar{X})$ (s'il existe).
\\

Pour le cas non-lin\'eaire, la condition (\ref{m9}) est n\'ecessaire mais non suffisante. Pour obtenir le minimum, il faut que la matrice $(\displaystyle \frac{\partial^2 G}{\partial X_i \partial X_j} ),i,j\in \left\{1,2,...,m\right\} $ soit d\'efinie positive.
\section{Etude d'un cas pratique}
On consid\`ere la d\'etermination d'un point par trilat\'eration dans l'option de la g\'eod\'esie tridimensionnelle pour la d\'etermination des coordonn\'ees $(x,y,z)$ d'un point $M$ inconnu \`a partir des mesures des distances vers $4$ points connus $M_k(u_k,v_k,w_k)_{k=1,2,3,4}$. 
\\

Pour faciliter les calculs, on prendra $P=I,\,\sigma_0=1$ et nous adoptons la fonction $\zeta(X)$ comme suit:
\begin{align}
	\zeta(X)=\left|
\begin{array}{l}
	\zeta_1=(x-u_1)^2+(y-v_1)^2+(z-w_1)^2 \\
		\zeta_2=(x-u_2)^2+(y-v_2)^2+(z-w_2)^2 \\
	\zeta_3=(x-u_3)^2+(y-v_3)^2+(z-w_3)^2 \\
		\zeta_4=(x-u_4)^2+(y-v_4)^2+(z-w_4)^2 
\end{array}\right. \label{m10} \\
\mbox{avec} \quad X_1=x,\quad X_2=y,\quad X_3=z \quad \mbox{les inconnues} \label{m11} 
\end{align}
et d'apr\`es (\ref{m10}), la fonction $\zeta$  n'est pas une fonction lin\'eaire des variables $X_i$. $\zeta$ est une fonction de $\BbR^3 \rightarrow \BbR^4$ qui s'\'ecrit:
$$ \zeta(X)=\zeta_ke_k$$ o\`u $e_k$ est la base orthonorm\'ee de $\BbR^n$. Voyons qu'elle v\'erifie les deux conditions (i) et (ii) cit\'es ci-dessus.
\\

Calculons  $\displaystyle \frac{\partial \zeta}{\partial X_i}$, on a alors:
\begin{equation}
	\frac{\partial \zeta}{\partial X_1}=\left|
\begin{array}{l}
	2(x-u_1) \\
	2(x-u_2) \\
	2(x-u_3) \\
	2(x-u_4)\\
\end{array}\right., \quad 	\frac{\partial \zeta}{\partial X_2}=\left|
\begin{array}{l}
	2(y-v_1) \\
	2(y-v_2) \\
	2(y-v_3) \\
	2(y-v_4)\\
\end{array}\right., \quad	\frac{\partial \zeta}{\partial X_3}=\left|
\begin{array}{l}
	2(z-w_1) \\
	2(z-w_2) \\
	2(z-w_3) \\
	2(z-w_4)\\
\end{array}\right. \quad \label{m12}
 \end{equation}
Pourque les 3 vecteurs soient lin\'eairement ind\'ependants, il faut que les points $M,M_i,M_j,M_k$ ne soient pas align\'es. Pour la condition (ii), on a facilement:
\begin{equation}
	\frac{\partial^2 \zeta}{\partial X^2_1}=\left|
\begin{array}{l}
	2 \\
	2 \\
	2 \\
	2
\end{array}\right., \quad 	\frac{\partial^2 \zeta}{\partial X^2_2}=\left|
\begin{array}{l}
	2 \\
	2 \\
	2 \\
	2
\end{array}\right., \quad	\frac{\partial^2 \zeta}{\partial X^2_3}=\left|
\begin{array}{l}
	2 \\
	2 \\
	2 \\
	2
\end{array}\right. \quad \label{m12a}
\end{equation}
et pour $i\neq j$, on a:
\begin{equation}
	\frac{\partial^2 \zeta}{\partial X_i \partial X_j}=\left|
\begin{array}{l}
	0 \\
	0 \\
	0 \\
	0
\end{array} \right. \label{m13}
\end{equation}
Donc les quantit\'es $\displaystyle \frac{\partial^2 \zeta(X)}{\partial X_i \partial X_j}$ sont continues et la condition (ii) est v\'erifi\'ee.
\subsection{Ecriture des Equations de Lagrange-Euler}
Pour d\'eterminer la solution par les moindres carr\'es du mod\`ele non-lin\'eaire, on \'ecrit les conditions (\ref{m8}). Le vecteur $L=(L_1,L_2,L_3,L_4)^T$ telque chacun des $L_i$ repr\'esente le carr\'e de la distance spatiale mesur\'ee. On a alors en utilisant (\ref{m12}):
$$ <L-\zeta(X), \frac{\partial \zeta(X)}{\partial X_i}>=0;\quad i=1,2,3 $$
soit:
\begin{align}
	(x-u_1)(L_1-\zeta_1)+(x-u_2)(L_2-\zeta_2)+(x-u_3)(L_3-\zeta_3)+(x-u_4)(L_4-\zeta_4)=0 \nonumber \\
	(y-v_1)(L_1-\zeta_1)+(y-v_2)(L_2-\zeta_2)+(y-v_3)(L_3-\zeta_3)+(y-v_4)(L_4-\zeta_4)=0 \nonumber \\
	(z-w_1)(L_1-\zeta_1)+(z-w_2)(L_2-\zeta_2)+(z-w_3)(L_3-\zeta_3)+(z-w_4)(L_4-\zeta_4)=0 \label{m14} 
	  	\end{align}
%\vspace{8mm}
Les \'equations (\ref{m14}) repr\'esente un syst\`eme de trois \'equations non lin\'eaires de trois inconnues $(x,y,z)$ dont la solution est un peu compliqu\'ee.
\\

Pour faciliter encore la r\'esolution du syst\`eme pr\'ec\'edent, on va supposer que la variable $z$ est connue \'egale \`a $z_0$, dans ce cas, on se limite \`a trois distances mesur\'ees $L_1, L_2$ et $L_2$. Alors (\ref{m14}) s'\'ecrit:
 \begin{align}
	(x-u_1)(L_1-\zeta_1)+(x-u_2)(L_2-\zeta_2)+(x-u_3)(L_3-\zeta_3)=0 \nonumber \\
	(y-v_1)(L_1-\zeta_1)+(y-v_2)(L_2-\zeta_2)+(y-v_3)(L_3-\zeta_3)=0 \label{m15} 
	  	\end{align}
avec:
$$ 	\zeta_1=(x-u_1)^2+(y-v_1)^2+(z_0-w_1)^2 $$
	$$	\zeta_2=(x-u_2)^2+(y-v_2)^2+(z_0-w_2)^2 $$
		$$	\zeta_3=(x-u_3)^2+(y-v_3)^2+(z_0-w_3)^2 $$
Les expressions $\zeta_i-L_i$ s'\'ecrivent sous la forme:
\begin{align}
			\zeta_i-L_i=x^2+y^2-2xu_i-2yv_i+a_i  \label{m16} \\
\mbox{avec} \quad a_i= constante  \label{m17} 
\end{align}
Le syst\`eme (\ref{m15}) devient:
\begin{align}
		(x-u_1)(x^2+y^2-2xu_1-2yv_1+a_1)+(x-u_2)(x^2+y^2-2xu_2-2yv_2+a_2) \nonumber \\
		+(x-u_3)(x^2+y^2-2xu_3-2yv_3+a_3)=0 \label{m18} \\
		(y-v_1)(x^2+y^2-2xu_1-2yv_1+a_1)+(y-v_2)(x^2+y^2-2xu_2-2yv_2+a_2) \nonumber \\
		+(y-v_3)(x^2+y^2-2xu_3-2yv_3+a_3)=0 \label{m19}
\end{align}
En d\'eveloppant les \'equations (\ref{m18}) et(\ref{m19}), on obtient:
\begin{align}
	3x^3+3xy^2-3x^2(u_1+u_2+u_3)-y^2(u_1+u_2+u_3)-2xy(v_1+v_2+v_3)+\nonumber \\ x(a_1+a_2+a_3+ 2u_1^2+2u_2^2+2u_3^2)+2y(u_1v_1+u_2v_2+u_3v_3)-(a_1u_1+a_2u_2+a_3u_3)=0 \label{m20} \\
		3y^3+3yx^2-3y^2(v_1+v_2+v_3)-x^2(v_1+v_2+v_3)-2xy(u_1+u_2+u_3)+\nonumber \\ y(a_1+a_2+a_3+ 2v_1^2+2v_2^2+2v_3^2)+2x(u_1v_1+u_2v_2+u_3v_3)-(a_1v_1+a_2v_2+a_3v_3)=0 \label{m21}
\end{align}
Supposons qu'on se limite \`a deux distances $L_1$ et $L_2$, alors on a \`a r\'esoudre :
	\begin{align}
	2x^3+2xy^2-2x^2(u_1+u_2)-y^2(u_1+u_2)-2xy(v_1+v_2)+\nonumber \\ x(a_1+a_2+ 2u_1^2+2u_2^2)+2y(u_1v_1+u_2v_2)-(a_1u_1+a_2u_2)=0 \label{m22} \\
		2y^3+2yx^2-2y^2(v_1+v_2)-x^2(v_1+v_2)-2xy(u_1+u_2)+\nonumber \\ y(a_1+a_2+ 2v_1^2+2v_2^2)+2x(u_1v_1+u_2v_2)-(a_1v_1+a_2v_2)=0 \label{m23}
\end{align}
\subsection{R\'eduction des Equations de Lagrange-Euler}
Dans ce paragraphe, on essaye de r\'eduire l'\'ecrirure du syst\`eme (\ref{m22}) - (\ref{m23}). A cet effet posons:
\begin{align}
	z=x+iy ,\quad \bar{z}=x-iy \quad \mbox{avec}\quad i=\sqrt{-1} \nonumber \\
	s=u_1+u_2,\quad t=v_1+v_2,\quad p=u_1^2+u_2^2,\quad q=v_1^2+v_2^2 \nonumber \\
	a=a_1+a_2,\quad r=u_1v_1+u_2v_2,\quad d=a_1u_1+a_2u_2,\quad f=a_1v_1+a_2v_2 \label{m24}
\end{align}
ce qui donne:
\begin{equation}
	2x=z+\bar{z},\quad 2iy=z-\bar{z}  \label{m25}
\end{equation}
Alors les expressions (\ref{m22}) - (\ref{m23}) deviennent:
\begin{align}
	4z^2\bar{z}+4z\bar{z}^2+(2it-s)z^2-6sz\bar{z}-(s+2it)\bar{z}^2 
	+2z(a+2p-2ir)+2\bar{z}(a+2p+2ir)-4d=0 \label{m26} \\
	4z^2\bar{z}-4z\bar{z}^2+(it-2s)z^2-6itz\bar{z}+(2s+it)\bar{z}^2 
	+2z(a+2q+2ir)-2\bar{z}(a+2q-2ir)-4if=0 \label{m27}
\end{align}
Nous pr\'esentons dans la suite la r\'esolution des \'equations ci-dessus:

\section{R\'esolution du Syst\`eme}
Posons:
\begin{align}
	z= x+iy=\rho e^{i \theta} \label{ma3} \\
	   s+2it=  l e^{i \omega} \label{ma4} \\
	a+2p+2ir=m   e^{i \alpha} \label{ma5} \\
	a+2q+2ir=k   e^{i \mu   } \label{ma6} \\
	2s+it  =h    e^{i \varphi}\label{ma7} 
\end{align}
Utilisant les relations:
\begin{align}
	z+\bar{z}=2\rho cos \theta \label{ma8} \\ 
	z- \bar{z}=2i\rho sin \theta ,\quad z\bar{z}=\rho^2 \label{ma9}
\end{align}
Les \'equations (\ref{m26}) et (\ref{m27}) deviennent:
\begin{align}
8\rho^3 cos\theta -6s\rho^2-2l\rho^2cos(2\theta - \omega)+4m\rho cos(\theta-\alpha)-4d=0 \label{ma10} \\
8i\rho^3 sin\theta -6it\rho^2-2ih\rho^2sin(2\theta - \varphi)+4ik\rho sin(\theta+\mu)-4if=0 \label{ma11} 
\end{align}
Soit:
\begin{align}
4\rho^3 cos\theta -\left[3s+lcos(2\theta - \omega)\right]\rho^2+2m\rho cos(\theta-\alpha)-2d=0 \label{ma12} \\
4\rho^3 sin\theta -\left[3t+hsin(2\theta - \varphi)\right]\rho^2+2k\rho sin(\theta+\mu)-2f=0 \label{ma13} 
\end{align}
Eliminant les termes constants entre (\ref{ma12}) et(\ref{ma13}) et apr\`es simplification par $\rho \neq 0$, on obtient l'\'equation:
\begin{align}
4\rho^2 (d sin \theta -f cos\theta)+\rho(3sf-3td+lfcos(2\theta-\omega)-hdsin(2\theta-\varphi))+\nonumber \\
2kdsin(\theta+\mu)-2mfcos(\theta-\alpha)=0 \label{ma14}
\end{align}
Maintenant, on \'elimine les coefficients en $\rho^3$ des \'equations (\ref{ma12}) et(\ref{ma13}), nous obtenons une \'equation en deuxi\`eme degr\'e en $\rho$:
\begin{align}
	\rho^2(3s.sin\theta -3tcos\theta+lsin\theta cos(2\theta-\omega)-hcos\theta sin(2\theta-\varphi))+\nonumber \\
	2\rho(kcos\theta sin(\theta+\mu)-msin\theta cos(\theta-\alpha))-2fcos\theta +2dsin\theta=0 \label{ma15}
\end{align}
\section{R\'esolution du syst\`eme des inconnues en $\rho$ et $\theta$}
On \'ecrit les \'equations (\ref{ma14}) et(\ref{ma15}) sous la forme:
\begin{align}
	A\rho^2+B\rho+C=0 \label{ma16} \\
	A'\rho^2+B'\rho+C'=0 \label{ma17}
\end{align}
Le syst\`eme pr\'ec\'edent est r\'esolvable si et seulement si:
\begin{equation}
	\frac{A}{A'}=\frac{B}{B'}=\frac{C}{C'} \label{ma18}
\end{equation}
soit en choisissant:
\begin{equation}
	\frac{B}{B'}=\frac{C}{C'} \Rightarrow \,\,BC'-B'C=0 \label{ma19}
\end{equation}
Soit:
\begin{align}
	\left[3sf-3td+lfcos(2\theta-\omega)-hdsin(2\theta-\varphi)\right](2dsin\theta-2fcos\theta)= \nonumber \\ \left[2kcos\theta sin(\theta+\mu)-2msin\theta cos(\theta-\alpha)\right]\left[2kdsin(\theta+\mu)-2mfcos(\theta-\alpha)\right] \label{ma20}
\end{align}
ce qui donne un polyn\^ome en $sin\theta, cos\theta, sin2\theta, cos2\theta$. En posant:
\begin{equation}
	\xi=tg\theta \label{ma21}
\end{equation}
et utilisant les formules:
\begin{align}
	sin2\theta=\frac{2\xi}{1+\xi^2} \label{ma22} \\
	cos2\theta=\frac{1-\xi^2}{1+\xi^2} \label{ma23}
\end{align}
L'\'equation (\ref{ma20}) devient un polyn\^ome du troisi\`eme degr\'e en $\xi$ qu'on peut r\'esoudre par les m\'ethodes classiques.
\section{Calcul de la Matrice Covariance des inconnues}
Comme on a:
\begin{equation}
	\zeta=\left(
\begin{array}{l}
 	\zeta_1=(x-u_1)^2+(y-v_1)^2+(z_0-w_1)^2 \\
		\zeta_2=(x-u_2)^2+(y-v_2)^2+(z_0-w_2)^2 \\
	\zeta_3=(x-u_3)^2+(y-v_3)^2+(z_0-w_3)^2 \\
\end{array} \right) \label{m28a}
\end{equation}
Or, on s'est limit\'e \`a deux inconnues, le vecteur $\zeta=(\zeta_1,\zeta_2)^T$. D'o\`u en diff\'erentiant (\ref{m28a}) en utilisant les deux premi\`eres lignes, on obtient alors:
\begin{equation}
	d\zeta=\left(
\begin{array}{l}
 	d\zeta_1                                 \\
	d	\zeta_2                                 
\end{array} \right) =\left(
\begin{array}{ll}
 	       2(x-u_1) & 2(y-v_1)               \\
		       2(x-u_2) & 2(y-v_2)               
\end{array} \right).\left(
\begin{array}{l}
	dx \\
	dy
\end{array}\right)=J_X.dX \label{m29a}
\end{equation}
avec:
\begin{equation}
	J_X=2.\left(
\begin{array}{ll}
	 	     x-u_1 & y-v_1               \\
		       x-u_2 & y-v_2               
\end{array} \right) \label{m30a}
\end{equation}
Par suite:
\begin{equation}
	dX=\left(
\begin{array}{l}
 	dx                                 \\
	dy                                 
\end{array} \right) =J_X^{-1}.\left(
\begin{array}{l}
 	       d\zeta_1  \\
		       d\zeta_2                
\end{array} \right) \label{m31a}
\end{equation}
d'o\`u la matrice covariance du vecteur inconnu  $X=(x,y)^T$:
\begin{equation}
	\sigma_X^2={J_X}^{-1}.\sigma_L^2.{{J_X}^{-1}}^T \label{m32a}
\end{equation}

\quad[2]\textbf{ P.J.G Teunissen}. The Geometry of Geodetic Inverse Linear Mapping and Non-Linear Adjustment. Publications On Geodesy, $n^{\circ}1$, Volume 8, Netherlands Geodetic Commission. 1985. 177p. 

\addtocontents{toc}{\protect\vspace{\baselineskip}}

%% The Bibliography
%% ==> You need a file 'literature.bib' for this.
%% ==> You need to run BibTeX for this (Project | Properties... | Uses BibTeX)
%\addcontentsline{toc}{chapter}{\textbf{R\'ef\'erences}\,\ldots\ldots\ldots\ldots\ldots\ldots\ldots\ldots\ldots\ldots\ldots\ldots\ldots\ldots\ldots\ldots\ldots\ldots\mbox{   }\,\,\,\,\,\,\,\,\,\,} %'Bibliography' into toc
%\nocite{*} %Even non-cited BibTeX-Entries will be shown.
\bibliographystyle{plain} %Style of Bibliography: plain / apalike / amsalpha / ...
%\bibliography{literature} %You need a file 'literature.bib' for this.

%% The List of Figures
\clearpage
%\addcontentsline{toc}{chapter}{Liste des Figures}
%\listoffigures

%% The List of Tables
\clearpage
%\addcontentsline{toc}{chapter}{Liste des Tables}
%\listoftables
%\begin{align}
%\end{align}

%%%%%%%%%%%%%%%%%%%%%%%%%%%%%%%%%%%%%%%%%%%%%%%%%%%%%%%%%%%%%
%% APPENDICES
%%%%%%%%%%%%%%%%%%%%%%%%%%%%%%%%%%%%%%%%%%%%%%%%%%%%%%%%%%%%%
%\appendix
%% ==> Write your text here or include other files.

%\input{FileName} %You need a file 'FileName.tex' for this.

\end{document}